\documentclass{birkjour}
\usepackage{hyperref, enumerate}
 \newtheorem{theorem}{Theorem}[section]
 \newtheorem{corollary}[theorem]{Corollary}
 \newtheorem{lemma}[theorem]{Lemma}
 \newtheorem{proposition}[theorem]{Proposition}
 \theoremstyle{definition}
 \newtheorem{definition}[theorem]{Definition}
 \theoremstyle{remark}
 \newtheorem{remark}[theorem]{Remark}
 \newtheorem{ex}[theorem]{Example}
 \numberwithin{equation}{section}
\begin{document}

%
%
%
%
%
\firstpage{1}
\issuenumber{3}
\received{November 08, 2022}
\revised{February 27, 2023}
\accepted{March 16, 2023}

%
%

\title[Conformal Submersions Whose Total Manifolds Admit a Ricci Soliton]
 {Conformal Submersions Whose Total\\ Manifolds Admit a Ricci Soliton}

\author[Kiran Meena]{Kiran Meena}
\address{Harish-Chandra Research Institute\\ A CI of Homi Bhabha National Institute\\ Chhatnag Road, Jhunsi, Prayagraj-211019, India.}
\email{kiranmeena@hri.res.in}
\author[Akhilesh Yadav]{Akhilesh Yadav}
\address{Department of Mathematics\\Institute of Science\\Banaras Hindu University\\Varanasi-221005, India.}
\email{akhilesha68@gmail.com}
\subjclass{Primary 53B20; Secondary 53C12; 53C25; 53C43}

\keywords{Riemannian manifold, Ricci soliton, Riemannian submersion, conformal submersion, Harmonic map.}


\begin{abstract}
In this paper, we study conformal submersions from Ricci solitons to Riemannian manifolds with non-trivial examples. First, we study some properties of the O'Neill tensor $A$ in the case of conformal submersion. We also find a necessary and sufficient condition for conformal submersion to be totally geodesic and calculate the Ricci tensor for the total manifold of such a map with different assumptions. Further, we consider a conformal submersion $F:M \to N$ from a Ricci soliton to a Riemannian manifold and obtain necessary conditions for the fibers of $F$ and the base manifold $N$ to be Ricci soliton, almost Ricci soliton and Einstein. Moreover, we find necessary conditions for a vector field and its horizontal lift to be conformal on $N$ and $(KerF_\ast)^\bot,$ respectively. Also, we calculate the scalar curvature of Ricci soliton $M$. Finally, we obtain a necessary and sufficient condition for $F$ to be harmonic.
\end{abstract}
\maketitle
\section{Introduction}\label{sec1}
Conformal submersions are a natural generalization of Riemannian submersions \cite{Neill_1966, Gunduzalp_2013, Gunduzalp_2019, Gunduzalp_2018}, which restricted to vectors orthogonal to its fibers is a horizontal conformal submersions (or horizontal conformal maps). The submersion is called Riemannian submersion if such restriction is an isometry \cite{Falcitelli_2008}, where the fundamental equations relating the curvatures of the total manifold and the base manifold of such submersion were given. In \cite{Gundmundsson_1992}, Gudmundsson obtained the fundamental equations for the conformal submersions. Further, horizontally conformal maps (conformal maps) were defined by Fuglede \cite{Fugledge_1978} and Ishihara \cite{Ishihara_1979} which are useful for characterization of harmonic morphisms and have applications in medical imaging (brain imaging) and computer graphics.

A smooth map $F:(M^m,g) \to (N^n,h)$ is said to be weakly conformal at $p\in M$ if there exists a number $\lambda^2(p)$ such that \cite{Baird_2003}
\begin{equation}\label{eqn1.1}
	h(F_\ast X, F_\ast Y) = \lambda^2(p) g(X, Y)~ \text{for}~ X, Y \in \Gamma(T_pM).
\end{equation} 
Let $(M^m,g)$ and $(N^n,h)$ be Riemannian manifolds then $F$ is called horizontally weakly conformal map at $p\in M$ if either $(i)$ $F_{\ast p} =0$, or $(ii)$ $F_{\ast p}$ maps the horizontal space $\mathcal{H}_p= (KerF_{\ast p})^\bot$ conformally onto $T_{F(p)} N$, i.e., $F_{\ast p}$ is surjective and satisfies (\ref{eqn1.1}) for $X, Y$ belong to $\mathcal{H}_p$. If a point $p$ is of type $(i)$, then it is called critical point of $F$ and if point $p$ is of type $(ii)$, then it is called regular point. The number $\lambda^2(p)$ is called the square dilation, it is necessarily non-negative and its square root $\lambda(p)$ is called the dilation. A horizontally weakly conformal map $F: M\to N$ is said to be horizontally homothetic if the gradient of its dilation $\lambda$ is vertical, i.e. $\mathcal{H} (grad~ \lambda)= 0$ at regular points. If a horizontally weakly conformal map $F$ has no critical points, then it is called horizontally conformal submersion \cite{Baird_2003}. Thus, a Riemannian submersion is a horizontally conformal submersion with dilation  identically one. 

The concept of harmonic maps and morphisms is a very useful tool for global analysis and differential geometry. The theory of harmonic maps has been developed in \cite{Eells_1964}, which is still an active field in differential geometry and it has applications to many different areas of mathematics and physics. A map between Riemannian manifolds is called harmonic if the divergence of its differential map vanishes. Since harmonic maps between Riemannian manifolds satisfy a system of quasi-linear partial differential equations, one would solve partial differential equations on certain manifolds. On the other hand, harmonic morphisms are maps between Riemannian manifolds which preserve germs of harmonic function, i.e. these (locally) pullback real valued harmonic functions to real valued harmonic functions. These are characterized as harmonic maps which are horizontally weakly conformal. Hence harmonic morphisms can be viewed as a subclass of harmonic maps \cite{Sahin_2010}. 

In 1988, the notion of Ricci soliton was introduced by Hamilton \cite{Hamilton_1988} to find a desired metric on a Riemannian manifold. For the metric on Riemannian manifold, the Ricci flow is an evolution equation (heat equation)
\begin{equation*}
	\frac{\partial}{\partial t} g(t) = -2 Ric.
\end{equation*}
Moreover, the self-similar solutions of Ricci flow are Ricci solitons which are natural generalization of an Einstein metric \cite{Besse_1987}. Let $(M, g)$ be a Riemannian manifold, if there exists a smooth vector field (potential vector field) $\xi$ which satisfies
\begin{equation}\label{eqn1.2}
	\frac{1}{2} (L_\xi g) + Ric + \mu g = 0,
\end{equation}
then $(M, g, \xi, \mu)$ is said to be Ricci soliton. Here $L_\xi g$ is the Lie derivative of the metric tensor $g$ with respect to $\xi$, $Ric$ is the Ricci tensor of $M$ and $\mu$ is a constant. Moreover, $\xi$ is called conformal vector field \cite{Deshmukh_2014} if it satisfies $L_{\xi} g = 2fg$, where $f$ is the potential function of $\xi$. A Ricci soliton $(M, g, \xi, \mu)$ is called shrinking, steady and expanding if $\mu <0,$ $\mu =0$ and $\mu>0$, respectively. In \cite{Pigola_2011}, Pigola et al. introduced almost Ricci soliton by taking $\mu$ as a variable function. In \cite{Perleman_2002}, Perelman used Ricci soliton to solve the Poincar\'e conjecture. Thereafter, the geometry of Ricci soliton has become the hot topic in the research due to the geometric importance and their wide applications in theoretical physics.

Recently, Riemannian submersions whose total spaces admitting Ricci soliton, $\eta$-Ricci soliton, almost $\eta$-Ricci-Bourguignon soliton, almost Yamabe soliton, $\eta$-Ricci-Yamabe soliton and conformal $\eta$-Ricci Soliton were studied in \cite{Meric_2019}, \cite{Gunduzalp_2020}, \cite{Bejan_2021}, \cite{Chaubev_2022}, \cite{Fatima_2021}, \cite{Meric_2020}, \cite{Siddiqi_2020} and \cite{Siddiqui_2022}. In addition, Riemannian maps whose total or base spaces admitting a Ricci soliton were studied in \cite{Yadav_Riem_total}, \cite{Yadav_Clairaut Riem_total}, \cite{GGupta_2022}, \cite{Yadav_Riem_base} and \cite{Yadav_Clairaut Riem_base}. 
In this paper, we study conformal submersions whose total manifolds admitting a Ricci soliton. The paper is organised as: in Sect. \ref{sec2}, we give some basic facts for the conformal submersion which are needed for this paper. In Sect. \ref{sec3}, we give some new results for conformal submersion. We obtain necessary and sufficient condition for conformal submersion to be totally geodesic. Moreover, we calculate the Ricci tensor for the total manifold of such a map. In Sect. \ref{sec4}, we obtain necessary conditions for the fibers of conformal submersion $F$ to be Ricci soliton, almost Ricci soliton and Einstein. Also, we obtain necessary conditions for the base manifold to be Ricci soliton, almost Ricci soliton and Einstein. Moreover, we calculate the scalar curvature of total manifold $M$ by using Ricci soliton. Finally, we discuss the harmonicity of conformal submersion whose total space is Ricci soliton. In Sect. \ref{sec5}, we give four non-trivial examples to support the theory of the paper.
\section{Preliminaries}\label{sec2}
In this section, we recall the notion of conformal submersion between Riemannian manifolds and give a brief review of basic facts.

A surjective smooth map $F:(M^m,g) \to (N^n,h)$ between Riemannian manifolds is said to be Riemannian submersion if it has maximal rank at every point $p\in M$. The fibers of $F$ over $q \in N$ is defined as $F^{-1}(q)$. The vectors tangent to fibers form the smooth vertical distribution denoted by $\nu_p$ and its orthogonal complementary with respect to $g$ is called horizontal distribution denoted by $\mathcal{H}_p$. Projections onto the horizontal and vertical distributions is denoted by $\mathcal{H}$ and $\nu$, respectively. A vector field $E$ on $M$ is said to be projectable if there exists a vector field $\tilde{E}$ on $N$ such that $F_{\ast p}(E) = \tilde{E}_{F(p)}$. Then $E$ and $\tilde{E}$ are called $F$-related. For all $\tilde{E}$ on $N$ there exists a unique vector field $E$ on $M$ such that $E$ and $\tilde{E}$ are $F$-related, and the vector field $E$ is called the horizontal lift of $\tilde{E}$. 

A Riemannian submersion $F$ is called a conformal (horizontally conformal) submersion if $F_\ast$ restricted to horizontal distribution of $F$ is a conformal map, that is there exists a smooth function $\lambda: M \to \mathbb{R}^{+}$ such that
\begin{equation}\label{eqn2.1}
	h(F_\ast X, F_\ast Y) = \lambda^2(p) g(X, Y),~ \forall X, Y \in \Gamma(KerF_\ast)^\bot~ \text{and}~ p \in M.
\end{equation}
The O'Neill tensors $A$ and $T$ defined in \cite{Neill_1966} as
\begin{equation}\label{eqn2.2}
	A_E E'= \mathcal{H} \nabla_{\mathcal{H}E} \nu E' + \nu \nabla_{\mathcal{H}E} \mathcal{H}E',
\end{equation}
\begin{equation}\label{eqn2.3}
	T_E E'= \mathcal{H} \nabla_{\nu E} \nu E' + \nu \nabla_{\nu E} \mathcal{H}E',
\end{equation}
$\forall E, E' \in \Gamma(TM)$, where $\nabla $ is the Levi-Civita connection of $g$. For any $ E \in \Gamma(TM) $, $T_E$ and $A_E$ are skew-symmetric operators on $(\Gamma(TM),g)$ reversing the horizontal and the vertical distributions. It is also easy to see that $ T$  is vertical, $T_E = T_{\nu E} $ and $A$ is horizontal, $A_E = A_{\mathcal{H}E} $. We note that the tensor field $T$ satisfies $T_U W = T_W U,~ \forall U,W \in \Gamma(KerF_\ast)$. Now, from (\ref{eqn2.2}) and (\ref{eqn2.3}), we have
\begin{equation}\label{eqn2.4}
	\nabla_U V = T_U V +  \nu \nabla_U V,
\end{equation}
\begin{equation}\label{eqn2.5}
	\nabla_X U = A_X U +  \nu \nabla_X U,
\end{equation}
\begin{equation}\label{eqn2.6}
	\nabla_X Y = A_X Y +  \mathcal{H} \nabla_X Y,
\end{equation}
$\forall X,Y \in \Gamma(KerF_\ast)^\bot$ and $ U, V \in \Gamma(KerF_\ast)$. A conformal submersion $F$ is with totally umbilical fibers if \cite{Zawadzki_2014, Zawadzki_2020}
\begin{equation}\label{eqn2.7}
	T_U V = g(U,V)H ~\text{or}~ T_U X= -g(H, X)U,
\end{equation}
$\forall U,V \in \Gamma(KerF_\ast)$ and $X \in \Gamma(KerF_\ast)^\bot$, where $H$ is the mean curvature vector field of the fibers.
\begin{proposition}\label{prop2.1}
	\cite{Gundmundsson_1992}
	Let $F: (M, g) \to (N, h)$ be a horizontally conformal submersion such that $(KerF_\ast)^\bot$ is integrable. Then the horizontal space is totally umbilical in $(M, g)$, i.e. $A_X Y = g(X,Y) H'~\forall X,Y \in \Gamma(KerF_\ast)^\bot$, where $H'$ is the mean curvature vector field of $(KerF_\ast)^\bot$ given by 
	\begin{equation}\label{eqn2.9}
		H'= - \frac{\lambda^2}{2} \left(\nabla_\nu \frac{1}{\lambda^2}\right).
	\end{equation}
\end{proposition}
\noindent The differential $F_\ast$ of $F$ can be viewed as a section of bundle $Hom(TM,F^{-1}TN)$ $\to M$, where $F^{-1}TN$ is the pullback bundle whose fibers at $p\in M$ is $(F^{-1}TN)_p = T_{F(p)}N$, $p\in M$. The bundle $Hom(TM,F^{-1}TN)$ has a connection $\nabla$ induced from the Levi-Civita connection ${\nabla}^M$ and the pullback connection  ${\nabla}^F$. Then the second fundamental form \cite{Nore_1966} of $F$ is given by $(\nabla F_\ast) (X,Y) = {\nabla}_X^F F_\ast Y - F_\ast({\nabla}_X^M Y),~\forall X,Y \in \Gamma(TM)$, or
\begin{equation}\label{eqn2.10}
	(\nabla F_\ast) (X,Y) = {\nabla}_{F_\ast X}^N F_\ast Y - F_\ast({\nabla}_X^M Y),~\forall X,Y \in \Gamma(TM).
\end{equation} 
Note that for the sake of simplicity we write $\nabla^M$ as $\nabla$.
\begin{lemma}\label{lem2.1}
	\cite{Gundmundsson_1992}
	Let $F: (M^m, g) \to (N^n, h)$ be a horizontally conformal submersion. Then
	\begin{equation*}
		\begin{array}{ll}
			F_\ast(\mathcal{H} \nabla_X Y) = \nabla_{F_\ast X}^N F_\ast Y + \frac{\lambda^2}{2} \left\{X (\frac{1}{\lambda^2})F_\ast Y + Y (\frac{1}{\lambda^2})F_\ast X - g(X, Y) F_\ast (grad_\mathcal{H} ~(\frac{1}{\lambda^2}))\right\},
		\end{array}
	\end{equation*}
	for $X, Y$ basic vector fields and $\nabla$ Levi-Civita connection  on $M$.
\end{lemma}
Now, we denote $R$, $R^N$ and $R^\nu$ the Riemannian curvature tensors of $M$, $N$ and fibers of $F$, respectively. Then for a horizontally conformal submersion, we have 	\cite{Gundmundsson_1992}
\begin{equation}\label{eqn2.12}
	g(R(U,V)W, S) = g(R^\nu (U, V)W, S) + g(T_U W, T_V S) - g(T_V W, T_U S),
\end{equation}
\begin{equation}\label{eqn2.13}
	g(R(U,V)W, X) = g((\nabla_U T)_V W, X)- g((\nabla_V T)_U W, X),
\end{equation}
\begin{eqnarray}\label{eqn2.14}
	\begin{array}{ll}
		g(R(U,X)Y, V) &= g((\nabla_U A)_X Y, V) + g(A_X U, A_Y V) - g((\nabla_X T)_U Y, V) \\&- g(T_V Y, T_U X) + \lambda^2 g(A_X Y, U) g(V, grad_\nu ~ (\frac{1}{\lambda^2})),
	\end{array}
\end{eqnarray}
\begin{equation}\label{eqn2.15}
	g(R(X, Y)Z, U) =  g((\nabla_X A)_Y Z, U) - g((\nabla_Y A)_X Z, U) - g(T_U Z, \nu [X,Y]),
\end{equation}
\begin{eqnarray}\label{eqn2.16}
	\begin{array}{ll}
		g(R(X, Y)Z, L) &= \frac{1}{\lambda^2} h(R^N(\tilde{X},\tilde{Y}) \tilde{Z}, \tilde{L}) + \frac{1}{4}\{  g(\nu [ X, Z], \nu [Y, L] ) \\&-g(\nu [ Y, Z], \nu [X, L]) + 2g(\nu [ X, Y], \nu [Z, L] ) \} \\& + \frac{\lambda^2}{2} \{ g(X, Z) g(\nabla_Y grad(\frac{1}{\lambda^2}), L)\\&- g(Y, Z) g(\nabla_X grad(\frac{1}{\lambda^2}), L)\\&+ g(Y, L) g(\nabla_X grad(\frac{1}{\lambda^2}), Z)\\&- g(X, L) g(\nabla_Y grad(\frac{1}{\lambda^2}), Z)\} \\& + \frac{\lambda^4}{4} \{( g(X, L) g(Y, Z) - g(Y, L) g(X, Z))\|grad(\frac{1}{\lambda^2})\|^2 \\& + g(X(\frac{1}{\lambda^2})Y - Y(\frac{1}{\lambda^2})X, L(\frac{1}{\lambda^2})Z - Z(\frac{1}{\lambda^2})L )\},
	\end{array} 
\end{eqnarray}
where $X,Y,Z,L \in \Gamma(KerF_\ast)^\bot$ and $U,V,W,S \in \Gamma(KerF_\ast)$. Also, $X, Y, Z$ and $L$ are the horizontal lift of $\tilde{X}, \tilde{Y}, \tilde{Z}$ and $\tilde{L}$, respectively. Moreover, $\nabla$ and $\nabla^N$ are Levi-Civita connections on $M$ and $N$, respectively. We denote the Ricci tensor and the scalar curvature by $Ric$ and $s$ respectively, defined as $Ric (X,Y)= trace (Z \mapsto R(Z,X)Y)$ and $s= trace Ric (X,Y)$ for $ X,Y \in \Gamma(TM)$.

Now we recall the gradient, divergence and Hessian \cite{Sahin_2017}. Let $f\in \mathcal{F}(M)$, then gradient of $f$, denoted by $\nabla f$ or $gradf$, given by
\begin{equation}\label{eqn2.17} 
	g (grad~f, X) = X(f),~for~X \in \Gamma(TM).
\end{equation}
Let $\{e_i\}_{1\leq i \leq m}$ be an orthonormal basis of $T_pM$ then
\begin{equation}\label{eqn2.18}
	g(X, Y) = \sum_{i=1}^{m} g(X, e_i) g(Y, e_i).
\end{equation}
The divergence of $X$, denoted by $div(X)$ and given by
\begin{equation}\label{eqn2.19}
	div(X)  = \sum_{i=1}^{m} g(\nabla_{e_i} X, e_i), ~\forall X\in \Gamma(TM).
\end{equation}
The Hessian tensor $h_f: \Gamma(TM) \to \Gamma(TM)$ of $f$ is given by
\begin{equation*}
	h_f(X) = \nabla_{X} \nabla f,~for~ X \in \Gamma(TM).
\end{equation*}
The Hessian form of $f$, denoted by $Hessf$ is given by
\begin{equation}\label{eqn2.21}
	Hessf (X,Y) = g(h_f(X), Y),~\forall X, Y \in \Gamma(TM).
\end{equation}
The Laplacian of $f \in \mathcal{F}(M)$, denoted by $\Delta f$, is given by
\begin{equation}\label{eqn2.22}
	\Delta f = div(\nabla f).	
\end{equation}
\begin{lemma}\label{lem2.2} 
	\cite{Gundmundsson_1992}
	Let $(M,g)$ be a Riemannian manifold and $f: M \to \mathbb{R}$ be a smooth function. Then
	\begin{equation*}
		g(\nabla_X grad(f), Y)= g(\nabla_Y grad(f), X), ~ \text{for}~ X, Y \in \Gamma(TM).
	\end{equation*}
\end{lemma}
\section{Characterizations of Conformal Submersion}\label{sec3}
In this section, we will find some interesting results for the conformal submersion which are useful to investigate its geometry.
\begin{proposition}\label{prop3.1}
	\cite{Gundmundsson_1992}
	Let $F: (M, g) \to (N, h)$ be a horizontally conformal submersion with dilation $\lambda$. Then
	\begin{equation*}
		A_X Y = \frac{1}{2} \left\{\nu [X, Y] - \lambda^2 g(X, Y) \left(\nabla_\nu \frac{1}{\lambda^2}\right)\right\}, ~\forall X, Y \in \Gamma(KerF_\ast)^\bot.
	\end{equation*}
\end{proposition}
\begin{theorem}\label{thm3.1}
	Let $F: (M, g) \to (N, h)$ be a horizontally conformal submersion with dilation $\lambda$ such that $(KerF_\ast)^\bot$ is totally geodesic. Then $\lambda$ is constant on $KerF_\ast$.
\end{theorem}
\begin{proof}
	By Proposition \ref{prop3.1}, we have
	\begin{equation}\label{eqn3.1}
		2 A_X Y = \nu [X, Y] - \lambda^2 g(X, Y) \left(\nabla_\nu \frac{1}{\lambda^2}\right),
	\end{equation}
	for all $X, Y \in \Gamma(KerF_\ast)^\bot$. Then by using (\ref{eqn3.1}), we get
	\begin{equation*}
		\begin{array}{ll}
			2 A_Y X = -\nu [X, Y] - \lambda^2 g(X, Y) (\nabla_\nu \frac{1}{\lambda^2}),
		\end{array}
	\end{equation*}
	which can be written as
	\begin{equation}\label{eqn3.2}
		2 A_Y X = -\left\{(\nu [X, Y] - \lambda^2 g(X, Y) \left(\nabla_\nu \frac{1}{\lambda^2}\right)\right\} -2\lambda^2 g(X, Y) \left(\nabla_\nu \frac{1}{\lambda^2}\right).
	\end{equation}
	Using (\ref{eqn3.1}) in (\ref{eqn3.2}),  we get
	\begin{equation}\label{eqn3.3}
		A_Y X = -A_X Y -\lambda^2 g(X, Y) \left(\nabla_\nu \frac{1}{\lambda^2}\right).
	\end{equation}
	Since $(KerF_\ast)^\bot$ is totally geodesic, (\ref{eqn3.3}) implies
	\begin{equation*}
		\nabla_\nu \frac{1}{\lambda^2}=0,
	\end{equation*}
	which completes the proof.
\end{proof}
\begin{remark}
	Note that for a conformal submersion $A$ is not alternating on the horizontal vector fields, while for a Riemannian submersion $A_X Y = -A_Y X,~ \forall X,Y \in \Gamma(KerF_\ast)^\bot$.
\end{remark}
\begin{theorem}\label{thm3.2}
	Let $F: (M, g) \to (N, h)$ be a horizontally conformal submersion. If $A$ is parallel then $\lambda$ is constant on $KerF_\ast$.
\end{theorem}
\begin{proof}
	For any $X \in \Gamma(KerF_\ast)^\bot$ and $W\in \Gamma(KerF_\ast)$, we have
	\begin{equation}\label{eqn3.4}
		g((\nabla_X A)_W X, W) = g(\nabla_X A_W X, W) - g(A_{\mathcal{H} \nabla_X W}X, W) - g(A_W A_X X, W).
	\end{equation}
	Since $A$ is horizontal tensor, $A_W =0$. Also, since $A$ is parallel, $\nabla_X A=0$. Then from (\ref{eqn3.4}), we get
	\begin{equation}\label{eqn3.5}
		-g(A_{\mathcal{H} \nabla_X W} X, W)=0.
	\end{equation}
	Using (\ref{eqn3.1}) in (\ref{eqn3.5}), we get
	\begin{equation*}
		g(A_X \mathcal{H} \nabla_X W, W)+ g\left(\lambda^2 g(\mathcal{H} \nabla_X W, X)\nabla_\nu \frac{1}{\lambda^2}, W\right) =0,
	\end{equation*}
	which implies
	\begin{equation*}
		g(A_X A_X W, W)+ \lambda^2 g(\mathcal{H} \nabla_X W, X)g\left(\nabla_\nu \frac{1}{\lambda^2}, W\right) = 0.
	\end{equation*}
	Since $A$ is parallel, above equation implies
	\begin{equation*}
		\left\{ g\left(\nabla_\nu \frac{1}{\lambda^2}, W\right)\right\}^2 =0.
	\end{equation*}
	This completes the proof.
\end{proof}

\begin{theorem}\label{thm3.3}
	Let $F: (M,g) \to (N,h)$ be a horizontally conformal submersion. Then $F$ is totally geodesic if and only if
	\begin{enumerate}[(i)]
		\item fibers of $F$ are totally geodesic,
		\item $(KerF_\ast)^\bot$ is totally geodesic,
		\item $F$ is homothetic.
	\end{enumerate}
\end{theorem}
\begin{proof}
	We know that a map $F: (M^m,g) \to (N^n,h)$ between Riemannian manifolds is totally geodesic $ \iff (\nabla F_\ast)(X, Y) =0, \forall X, Y \in \Gamma(TM) \iff$ 
	
	\noindent $(\nabla F_\ast)(U, V) =0,~ \forall U, V \in \Gamma(KerF_\ast)$, $(\nabla F_\ast)(X, U) =0,~ \forall U \in \Gamma(KerF_\ast)$ and $X \in \Gamma(KerF_\ast)^\bot$, and
	$(\nabla F_\ast)(X, Y) =0,~ \forall X, Y \in \Gamma(KerF_\ast)^\bot$.
	Since $(\nabla F_\ast)(U, V) \in \Gamma(RangeF_\ast) = \Gamma(TN)$, $(\nabla F_\ast)(U, V) =0 \iff h((\nabla F_\ast)(U,$ $V), F_\ast Y) =0 ~\text{for} ~Y \in \Gamma(KerF_\ast)^\bot \iff h((\nabla_{F_\ast U}^N F_\ast V), F_\ast Y) =0 \iff h(F_\ast(\nabla_U V),$ $F_\ast Y)=0 \iff h(F_\ast(\mathcal{H} \nabla_U V), F_\ast Y)=0 \iff \lambda^2 g(\mathcal{H} \nabla_U V, Y) =0 \iff g(T_U V, Y)=0~ \forall Y \in \Gamma(KerF_\ast)^\bot \iff T_U V =0 \iff$ fibers of $F$ are totally geodesic. 
	
	Further, since $(\nabla F_\ast)(X, U) \in \Gamma(RangeF_\ast) = \Gamma(TN)$, $(\nabla F_\ast)(X, U) =0 \iff h((\nabla F_\ast)(X, U), F_\ast Y) =0 ~\text{for} ~Y \in \Gamma(KerF_\ast)^\bot \iff h((\nabla_{F_\ast X}^N F_\ast U), F_\ast Y) =0 \iff h(F_\ast(\nabla_X U), F_\ast Y)=0 \iff h(F_\ast(\mathcal{H} \nabla_X U), F_\ast Y)=0 \iff \lambda^2 g(\mathcal{H} \nabla_X U, Y) =0 \iff g(A_X U, Y)=0 \iff g(A_X Y, U)=0~ \forall U \in \Gamma(KerF_\ast) \iff A_X Y =0 \iff (KerF_\ast)^\bot$ is totally geodesic.
	
	Finally, since $(\nabla F_\ast)(X, Y) \in \Gamma(RangeF_\ast) = \Gamma(TN)$, $(\nabla F_\ast)(X, Y) =0 \iff h((\nabla F_\ast)(X, Y), F_\ast Z) =0 ~\text{for} ~Z \in \Gamma(KerF_\ast)^\bot \iff h((X~ \ln \lambda) F_\ast Y + (Y~ \ln \lambda) F_\ast X - g(X,Y) F_\ast (\nabla_{\mathcal{H}} \ln \lambda), F_\ast Z) =0 \iff \lambda^2 \{(X~\ln \lambda) g(Y, Z) + (Y~ \ln \lambda) g(X, Z) - g(X, Y)$ $g(\nabla_{\mathcal{H}} \ln \lambda, Z)\} =0 \iff (X~\ln \lambda) g(Y, Z) + (Y~ \ln \lambda) g(X, Z) - g(X, Y) g(\nabla_{\mathcal{H}} \ln \lambda, Z) =0 \iff F$ is homothetic.	 
\end{proof}
\begin{lemma}\label{lem3.1}
	Let $F: (M^m,g) \to (N^n,h)$ be a horizontally conformal submersion with dilation $\lambda$ such that $(KerF_\ast)^\bot$ is integrable. Then following statements are true:
	\begin{enumerate}[(i)]
		\item $g(A_{X_j}U, A_{X_j} V)= n^2 \frac{\lambda^4}{4} g(\nabla_\nu \frac{1}{\lambda^2}, U)g(\nabla_\nu \frac{1}{\lambda^2}, V)$,
		\item $g((\nabla_U A)_{X_j} X_j, V) = n g(\nabla_U H', V)$,
		\item $g((\nabla_X A)_{X_j} X_j, U) = n g(\nabla_X H', U)$,
		\item $g((\nabla_{X_j} A)_X X_j, U)= g(X, X_j) g(\nabla_{X_j} H', U)$,
		\item $g((\nabla_{U_i} A)_X Y, U_i)= g(X, Y) ~div( H')$,		
		\item $g(A_X U_i, A_Y U_i)= g(X, Y) \frac{\lambda^4}{4} |\nabla_\nu \frac{1}{\lambda^2}|^2$,
	\end{enumerate}
	\noindent for $U, V \in \Gamma(KerF_\ast)~ and ~X, Y \in \Gamma(KerF_\ast)^\bot$, where $\{X_j\}_{1 \leq j \leq n}$ and $\{U_i\}_{n+1 \leq i \leq m}$ are orthonormal bases of $(Ker F_\ast)^\bot$ and $KerF_\ast$, respectively.
\end{lemma}
\begin{proof}
	Since $A_{X_j} U \in \Gamma(Ker F_\ast)^\bot$ for $U \in \Gamma(KerF_\ast)$. Then we can write
	\begin{equation*}
		A_{X_j} U= \sum\limits_{j=1}^{n} g(A_{X_j} U, X_j) X_j=- \sum\limits_{j=1}^{n} g(A_{X_j} X_j, U) X_j,
	\end{equation*}
	where $\{X_j\}_{1 \leq j \leq n}$ is an orthonormal basis of $(Ker F_\ast)^\bot$. Then by using above equation, we get
	\begin{equation}\label{eqn3.7}
		g(A_{X_j} U, A_{X_j} V)= \sum\limits_{j=1}^{n} g(A_{X_j} X_j, U)g(A_{X_j} X_j, V).
	\end{equation}
	Since $(KerF_\ast)^\bot$ is integrable, using Proposition \ref{prop2.1} in (\ref{eqn3.7}), we get 
	\begin{equation}\label{eqn3.8}
		g(A_{X_j} U, A_{X_j} V)= n^2 g(H', U) g(H', V).
	\end{equation}
	Using (\ref{eqn2.9}) in (\ref{eqn3.8}), we get the proof of the first statement.
	
	\noindent Also, since
	\begin{equation}\label{eqn3.9}
		g((\nabla_U A)_{X_j} X_j, V) = g(\nabla_U A_{X_j} X_j - A(\nabla_U X_j, X_j) - A(X_j, \nabla_U X_j), V).
	\end{equation}
	Since $(KerF_\ast)^\bot$ is integrable, using Proposition \ref{prop2.1} in (\ref{eqn3.9}), we get 
	\begin{equation}\label{eqn3.10}
		\begin{array}{ll}
			g((\nabla_U A)_{X_j} X_j, V) &= g(\nabla_U g(X_j, X_j) H' - g(\mathcal{H} \nabla_U X_j, X_j) H' \\&- g(X_j, \mathcal{H} \nabla_U X_j) H', V).
		\end{array}
	\end{equation}
	Using metric compatibility condition in (\ref{eqn3.10}), we get the proof of the second statement.
	\noindent Similarly, we can get the proof of the third, fourth, fifth and sixth statements easily.
\end{proof}
\begin{proposition}\label{prop3.2}
	Let $F: (M^m,g) \to (N^n,h)$ be a horizontally conformal submersion with dilation $\lambda$. Then
	\begin{eqnarray}
		\begin{array}{ll}\label{eqn3.11}
			Ric(U, V)&= Ric^\nu (U, V) - (m-n)g(T_U V, H) + \sum\limits_{j=1}^{n} g((\nabla_U A)_{X_j} X_j, V) \\&+ \sum\limits_{j=1}^{n} g(A_{X_j} U, A_{X_j}V) - \sum\limits_{j=1}^{n} g((\nabla_{X_j}T)_U X_j, V) \\&-\frac{\lambda^4}{2} n g(U, \nabla_\nu \frac{1}{\lambda^2})g(V, \nabla_\nu \frac{1}{\lambda^2}),
		\end{array}
	\end{eqnarray}
	\begin{equation}\label{eqn3.12}
		\begin{array}{ll}
			Ric(U, X)&= (m-n)g(\nabla_U H, X) -\sum\limits_{i=n+1}^{m} g((\nabla_{U_i} T)_U U_i, X) \\&+ \sum\limits_{j=1}^{n} g((\nabla_X A)_{X_j} X_j, U) - \sum\limits_{j=1}^{n} g((\nabla_{X_j} A)_X X_j, U) \\&- \sum\limits_{j=1}^{n} g(T_U X_j, \nu [X, X_j]),
		\end{array}
	\end{equation}
	\begin{eqnarray}\label{eqn3.13}
		\begin{array}{ll}
			Ric(X, Y)&= \sum\limits_{i=n+1}^{m} g((\nabla_{U_i} A)_X Y, U_i) + \sum\limits_{i=n+1}^{m} g(A_X U_i, A_Y U_i) \\&- \sum\limits_{i=n+1}^{m} g((\nabla_X T)_{U_i} Y, U_i) - \sum\limits_{i=n+1}^{m} g(T_{U_i} X, T_{U_i} Y) \\&+ \lambda^2 g(A_X Y, \nabla_\nu \frac{1}{\lambda^2}) + \frac{1}{\lambda^2} Ric^N (\tilde{X}, \tilde{Y}) \\&+\frac{3}{4} \sum\limits_{j=1}^{n} g(\nu [X, X_j], \nu [X_j, Y]) - \frac{(n-2)}{2} \lambda^2 g(\nabla_X \nabla \frac{1}{\lambda^2}, Y) \\& - \frac{\lambda^2}{2} g(X, Y) \left\{ \Delta^\mathcal{H}\frac{1}{\lambda^2} - n \left( H' \frac{1}{\lambda^2}\right)\right\} + \frac{ n \lambda^4}{4} g(X, Y) |\nabla \frac{1}{\lambda^2}|^2 \\&+ \frac{\lambda^4}{4} (n-2) (X \frac{1}{\lambda^2}) (Y\frac{1}{\lambda^2}),
		\end{array}
	\end{eqnarray}
	for $U, V \in \Gamma(KerF_\ast)~ and ~X, Y \in \Gamma(KerF_\ast)^\bot$, where $\{X_j\}_{1 \leq j \leq n}$ and $\{U_i\}_{n+1 \leq i \leq m}$ are orthonormal bases of $(Ker F_\ast)^\bot$ and $KerF_\ast$, respectively. Also, $X$ and $Y$ are the horizontal lift of $\tilde{X}$ and $\tilde{Y}$, respectively.
\end{proposition}
\begin{proof} Since we know that
	\begin{equation*}
		Ric(U, V) = \sum\limits_{i=n+1}^{m} g(R(U_i, U)V, U_i) + \sum\limits_{j=1}^{n} g(R(X_j, U)V, X_j), ~\text{for}~ U, V \in \Gamma(KerF_\ast),
	\end{equation*}
	where $\{X_j\}_{1 \leq j \leq n}$ and $\{U_i\}_{n+1 \leq i \leq m}$ are orthonormal bases of $(Ker F_\ast)^\bot$ and $KerF_\ast$, respectively. Then by using (\ref{eqn2.12}) and (\ref{eqn2.14}) in above equation, we get
	\begin{eqnarray}\label{eqn3.14}
		\begin{array}{ll}
			Ric(U, V) &  =\sum\limits_{i=n+1}^{m} g(R^\nu (U_i, U)V, U_i) + \sum\limits_{i=n+1}^{m} g(T_{U_i} V, T_U U_i) \\&- \sum\limits_{i=n+1}^{m}g(T_U V, T_{U_i} U_i) + \sum\limits_{j=1}^{n} g((\nabla_U A)_{X_j} X_j, V) \\&+ \sum\limits_{j=1}^{n} g(A_{X_j} U, A_{X_j} V) - \sum\limits_{j=1}^{n} g((\nabla_{X_j} T)_U X_j, V) \\& - \sum\limits_{j=1}^{n} g(T_V X_j, T_U X_j) + \lambda^2 \sum\limits_{j=1}^{n}g (A_{X_j} X_j, U) g(V, \nabla_\nu \frac{1}{\lambda^2}).
		\end{array}
	\end{eqnarray}
	Using (\ref{eqn3.3}) and $g(T_{U_i} V, T_U U_i) = g(T_V X_j, T_U X_j)$ in (\ref{eqn3.14}), we get (\ref{eqn3.11}).
	
	\noindent Also, since we know that
	\begin{equation*}
		Ric(U, X) = \sum\limits_{i=n+1}^{m}g(R(U_i, U)X, U_i) + \sum\limits_{j=1}^{n} g(R(X_j, U)X, X_j),
	\end{equation*}
	for $U \in \Gamma(KerF_\ast)$ and $X \in \Gamma(KerF_\ast)^\bot$, where $\{X_j\}_{1 \leq j \leq n}$ and $\{U_i\}_{n+1 \leq i \leq m}$ are orthonormal bases of $(Ker F_\ast)^\bot$ and $KerF_\ast$, respectively. Then by using (\ref{eqn2.13}) and (\ref{eqn2.15}) in above equation, we get (\ref{eqn3.12}). Finally, since
	\begin{equation*}
		Ric(X, Y) = \sum\limits_{i=n+1}^{m}g(R(U_i, X)Y, U_i) + \sum\limits_{j=1}^{n} g(R(X_j, X)Y, X_j),~\text{for}~X, Y \in \Gamma(KerF_\ast)^\bot,
	\end{equation*}
	where $\{X_j\}_{1 \leq j \leq n}$ and $\{U_i\}_{n+1 \leq i \leq m}$ are orthonormal bases of $(Ker F_\ast)^\bot$ and $KerF_\ast$, respectively. By using (\ref{eqn2.14}) and (\ref{eqn2.16}) in above equation, then using (\ref{eqn2.18}), (\ref{eqn2.17}) and Lemma \ref{lem2.2}, we get
	\begin{eqnarray}\label{eqn3.15}
		\begin{array}{ll}
			Ric(X, Y)&= \sum\limits_{i=n+1}^{m} g((\nabla_{U_i} A)_X Y, U_i) + \sum\limits_{i=n+1}^{m} g(A_X U_i, A_Y U_i) \\&- \sum\limits_{i=n+1}^{m} g((\nabla_X T)_{U_i} Y, U_i) - \sum\limits_{i=n+1}^{m} g(T_{U_i} X, T_{U_i} Y) \\&+ \lambda^2 g(A_X Y, \nabla_\nu \frac{1}{\lambda^2}) + \frac{1}{\lambda^2}\sum\limits_{j=1}^{n} h(R^N(\tilde{X_j}, \tilde{X})\tilde{Y}, \tilde{X_j}) \\&+\frac{3}{4} \sum\limits_{j=1}^{n} g(\nu [X, X_j], \nu [X_j, Y]) - (\frac{n}{2}-1) \lambda^2 g(\nabla_X \nabla \frac{1}{\lambda^2}, Y) \\& - \frac{\lambda^2}{2} \sum\limits_{j=1}^{n} g(X, Y) g(\nabla_{X_j} \nabla \frac{1}{\lambda^2}, X_j) + \frac{ (n-1) \lambda^4}{4} g(X, Y) |\nabla \frac{1}{\lambda^2}|^2 \\&+ \frac{\lambda^4}{4} \left\{ g(X, Y) |\nabla \frac{1}{\lambda^2}|^2 -2 (X \frac{1}{\lambda^2}) (Y\frac{1}{\lambda^2}) + n (X \frac{1}{\lambda^2}) (Y\frac{1}{\lambda^2}) \right\},
		\end{array}
	\end{eqnarray}
	where $\{\tilde{X_j}\}_{1\leq j\leq n}$ is an orthonormal basis of $TN$. Then using (\ref{eqn2.21}), (\ref{eqn2.22}) and (\ref{eqn2.9}) in (\ref{eqn3.15}), we get (\ref{eqn3.13}).
\end{proof}
\begin{corollary}\label{Cor3.1}
	Let $F: (M^m,g) \to (N^n,h)$ be a horizontally conformal submersion with dilation $\lambda$ such that fibers of $F$ are totally geodesic and $(Ker F_\ast)^\bot$ is integrable. Then
	\begin{equation*}
		\begin{array}{ll}
			Ric(U, V) = Ric^\nu (U, V) +ng(\nabla_U H', V) + (\frac{n^2}{4} - \frac{n}{2})\lambda^4 g(U, \nabla_\nu \frac{1}{\lambda^2})g(V, \nabla_\nu \frac{1}{\lambda^2}),
		\end{array}
	\end{equation*}
	\begin{equation*}
		Ric(U, X)= ng(\nabla_X H', U) - \sum\limits_{j=1}^{n} g(X, X_j) g(\nabla_{X_j} H', U),
	\end{equation*}
	\begin{equation*}
		\begin{array}{ll}
			Ric(X, Y)=& g(X, Y) div(H')+ \frac{1}{\lambda^2} Ric^N (\tilde{X}, \tilde{Y}) -\frac{3}{4} \lambda^4 g(X, Y) |\nabla_{\nu} \frac{1}{\lambda^2}|^2 \\&- \frac{(n-2)}{2} \lambda^2 g(\nabla_X \nabla \frac{1}{\lambda^2}, Y) - \frac{\lambda^2}{2} g(X, Y) \left\{ \Delta^\mathcal{H}\frac{1}{\lambda^2} - n \left( H' \frac{1}{\lambda^2}\right)\right\} \\&+ \frac{ n \lambda^4}{4} g(X, Y) |\nabla \frac{1}{\lambda^2}|^2 + \frac{\lambda^4}{4} (n-2) (X \frac{1}{\lambda^2}) (Y\frac{1}{\lambda^2}).
		\end{array}
	\end{equation*}
\end{corollary}
\begin{proof}
	Using Lemma \ref{lem3.1} and Proposition \ref{prop3.2}, we get the proof.
\end{proof}    
\begin{corollary}\label{cor3.2}
	Let $F: (M^m,g) \to (N^n,h)$ be a horizontally homothetic conformal submersion with dilation $\lambda$ such that the fibers of $F$ are totally geodesic and $(Ker F_\ast)^\bot$ is integrable. Then
	\begin{equation*}
		\begin{array}{ll}
			Ric(X, Y)=& g(X, Y) div(H')+ \frac{1}{\lambda^2} Ric^N (\tilde{X}, \tilde{Y}) \\&- \frac{1}{4} \lambda^4 g(X, Y) |\nabla_{\nu} \frac{1}{\lambda^2}|^2 + \frac{n \lambda^2}{2} g(X, Y) \left(H' \frac{1}{\lambda^2}\right).
		\end{array}
	\end{equation*}
\end{corollary}
\begin{proof}
	Using Lemma \ref{lem3.1}, Proposition \ref{prop3.2} and $\mathcal{H} (grad~ \lambda)= 0$, we get the proof.
\end{proof}
\begin{corollary}\label{cor3.3}
	Let $F: (M^m,g) \to (N^n,h)$ be a totally geodesic horizontally conformal submersion with dilation $\lambda$. Then
	\begin{equation*}
		Ric(U, V) = Ric^\nu (U, V),
	\end{equation*}
	\begin{equation*}
		Ric(U, X)=0,
	\end{equation*}
	\begin{equation*}
		\begin{array}{ll}
			Ric(X, Y)= \frac{1}{\lambda^2} Ric^N (\tilde{X}, \tilde{Y}).
		\end{array}
	\end{equation*}
\end{corollary}
\begin{proof}
	Using Theorems \ref{thm3.3} and \ref{thm3.1}, and Proposition \ref{prop3.2}, we get the proof.
\end{proof}
\begin{theorem}\label{thm3.4}
	Let $F: (M^m,g) \to (N^n,h)$ be a totally geodesic horizontally conformal submersion with dilation $\lambda$. Then
	\begin{equation*}
		s= s^{KerF_\ast}+ \frac{1}{\lambda^2}s^N,
	\end{equation*}
	where $s$, $s^{KerF_\ast}$ and $s^N$ denote the scalar curvatures of $M$, $KerF_\ast$ and $N$, respectively.
\end{theorem}
\begin{proof}
	Since the scalar curvature of $M$ is defined as \cite{Falcitelli_2008}
	\begin{equation}\label{eqn3.16}
		s = \sum\limits_{i=n+1}^{m} Ric(U_i, V_i) + \sum\limits_{j=1}^{n} Ric(X_j, X_j),
	\end{equation}
	where $\{X_j\}_{1 \leq j \leq n}$ and $\{U_i\}_{n+1 \leq i \leq m}$ are orthonormal bases of $(Ker F_\ast)^\bot$ and $KerF_\ast$, respectively. Then using Corollary \ref{cor3.3} in (\ref{eqn3.16}), we get the proof.
\end{proof}
\section{Conformal Submersion from Ricci Soliton}\label{sec4}
In this section, we consider a conformal submersion from a Ricci soliton to a Riemannian manifold and investigate its geometry.
\begin{theorem}
	Let $(M, g, \xi, \mu)$ be a Ricci soliton with the potential vector field $\xi \in \Gamma(TM)$ and $F:(M^m, g) \to (N^n, h)$ be a horizontally conformal submersion between Riemannian manifolds such that $KerF_\ast$ and $(KerF_\ast)^\bot$ are totally geodesic. Then following statements are true:
	\begin{enumerate}[(i)]
		\item If $\xi= W \in \Gamma(KerF_\ast)$, then any fiber of $F$ is a Ricci soliton.
		\item If $\xi= X \in \Gamma(KerF_\ast)^\bot,$ then any fiber of $F$ is an Einstein.
	\end{enumerate}
\end{theorem}
\begin{proof}
	Since $(M, g, \xi, \mu)$ is a Ricci soliton, then by (\ref{eqn1.2}), we have
	\begin{equation*}
		\frac{1}{2} (L_\xi g)(U,V) + Ric(U,V) + \mu g(U,V) = 0,~\text{for}~ U, V \in \Gamma(KerF_\ast),
	\end{equation*}
	which can be written as
	\begin{equation}\label{eqn4.1}
		\frac{1}{2} \{g(\nabla_U \xi, V) + g(\nabla_V \xi, U)\} + Ric(U,V) + \mu g(U,V) = 0.
	\end{equation}
	Using (\ref{eqn3.11}) and (\ref{eqn3.3}) in (\ref{eqn4.1}), we get
	\begin{eqnarray}\label{eqn4.2}
		\begin{array}{ll}
			\frac{1}{2} \{g(\nabla_U \xi, V) + g(\nabla_V \xi, U)\} + Ric^\nu (U, V) - (m-n)g(T_U V, H) \\+ \sum\limits_{j=1}^{n} g((\nabla_U A)_{X_j} X_j, V) + \sum\limits_{j=1}^{n} g(A_{X_j} U, A_{X_j}V) - \sum\limits_{j=1}^{n} g((\nabla_{X_j}T)_U X_j, V) \\+\lambda^2 g(A_{X_j} X_j, U)g(V, \nabla_\nu \frac{1}{\lambda^2}) + \mu g(U,V) = 0.
		\end{array}
	\end{eqnarray}
	Since $KerF_\ast$ and $(KerF_\ast)^\bot$ are totally geodesic, from (\ref{eqn4.2}), we get
	\begin{equation}\label{eqn4.3}
		\frac{1}{2} \{g(\nabla_U \xi, V) + g(\nabla_V \xi, U)\} + Ric^\nu (U, V)  + \mu g(U,V) = 0.
	\end{equation}
	If $\xi= W \in \Gamma(KerF_\ast)$, then from (\ref{eqn4.3}), we get 
	\begin{equation*}
		\frac{1}{2} \{g(\nabla_U W, V) + g(\nabla_V W, U)\} + Ric^\nu (U, V)  + \mu g(U,V) = 0,
	\end{equation*}
	which implies any fiber of $F$ is a Ricci soliton, this proves $(i)$.\\Now, if $\xi= X \in \Gamma(KerF_\ast)^\bot,$ then (\ref{eqn4.3}) implies 
	\begin{equation*}
		\frac{1}{2} \{g(\nabla_U X, V) + g(\nabla_V X, U)\} + Ric^\nu (U, V)  + \mu g(U,V) = 0.
	\end{equation*}
	Using (\ref{eqn2.4}) in above equation, we get
	\begin{equation}\label{eqn4.4}
		-g(T_U V, X) + Ric^\nu (U, V)  + \mu g(U,V) = 0.
	\end{equation}
	Since $KerF_\ast$ is totally geodesic, from (\ref{eqn4.4}), we get
	\begin{equation*}
		Ric^\nu (U, V)  + \mu g(U,V) = 0,
	\end{equation*}
	which implies $(ii)$. This completes the proof.
\end{proof}
\begin{theorem}
	Let $(M, g, \xi, \mu)$ be a Ricci soliton with the potential vector field $\xi \in \Gamma(TM)$ and $F:(M^m, g) \to (N^n, h)$ be a horizontally conformal submersion between Riemannian manifolds such that fibers of $F$ are totally umbilical and $(KerF_\ast)^\bot$ is totally geodesic. Then following statements are true:
	\begin{enumerate}[(i)]
		\item If $\xi= W \in \Gamma(KerF_\ast)$, then any fiber of $F$ is an almost Ricci soliton.
		\item If $\xi= X \in \Gamma(KerF_\ast)^\bot,$ then any fiber of $F$ is an Einstein.
	\end{enumerate}
\end{theorem}
\begin{proof}
	Since $(M, g, \xi, \mu)$ is a Ricci soliton and $(KerF_\ast)^\bot$ is totally geodesic, then by (\ref{eqn4.2}), we get
	\begin{equation*}
		\begin{array}{ll}
			\frac{1}{2} \{g(\nabla_U \xi, V) + g(\nabla_V \xi, U)\} + Ric^\nu (U, V) - (m-n)g(T_U V, H) \\- \sum\limits_{j=1}^{n} g((\nabla_{X_j}T)_U X_j, V) + \mu g(U,V) = 0.
		\end{array}
	\end{equation*}
	Since fibers of $F$ are totally umbilical, using (\ref{eqn2.7}) in above equation, we get
	\begin{equation}
		\begin{array}{ll}\label{eqn4.5}
			\frac{1}{2} \{g(\nabla_U \xi, V) + g(\nabla_V \xi, U)\} + Ric^\nu (U, V) - (m-n)g(U, V)\|H\|^2 \\+ \sum\limits_{j=1}^{n} g(\nabla_{X_j}H, X_j) g(U, V) + \mu g(U,V) = 0.
		\end{array}
	\end{equation}
	Using (\ref{eqn2.19}) in (\ref{eqn4.5}), we get
	\begin{equation}\label{eqn4.6}
		\begin{array}{ll}
			\frac{1}{2} \{g(\nabla_U \xi, V) + g(\nabla_V \xi, U)\} + Ric^\nu (U, V) - (m-n)\|H\|^2g(U, V) \\+ div(H) g(U, V) + \mu g(U,V) = 0.
		\end{array}
	\end{equation}
	If $\xi= W \in \Gamma(KerF_\ast)$, then (\ref{eqn4.6}) implies
	\begin{equation*}
		\frac{1}{2} \{g(\nabla_U W, V) + g(\nabla_V W, U)\} + Ric^\nu (U, V)+ f_1 g(U, V) = 0,
	\end{equation*}
	where $f_1 = div(H)- (m-n)\|H\|^2 + \mu$ is a smooth function on $M$. Thus any fiber of $F$ is an almost Ricci soliton, which implies $(i)$.\\Now, if $\xi= X \in \Gamma(KerF_\ast)^\bot,$ then (\ref{eqn4.6}) becomes 
	\begin{equation*}
		\begin{array}{ll}
			\frac{1}{2} \{g(\nabla_U X, V) + g(\nabla_V X, U)\} + Ric^\nu (U, V) - (m-n)\|H\|^2g(U, V) \\+ div(H) g(U, V) + \mu g(U,V) = 0.
		\end{array}
	\end{equation*}
	Using metric compatibility and (\ref{eqn2.4}) in above equation, we get
	\begin{equation*}
		\begin{array}{ll}
			-g(T_U V, X) + Ric^\nu (U, V) - (m-n)\|H\|^2g(U, V) + div(H) g(U, V) + \mu g(U,V) = 0.
		\end{array}
	\end{equation*}
	Since fibers of $F$ are totally umbilical, using (\ref{eqn2.7}) in above equation, we get
	\begin{equation*}
		\begin{array}{ll}
			-g(U, V)g(H, X) + Ric^\nu (U, V) - (m-n)\|H\|^2g(U, V) \\+ div(H) g(U, V) + \mu g(U,V) = 0,
		\end{array}
	\end{equation*}
	which implies
	\begin{equation*}
		Ric^\nu (U, V) + f_2 g(U, V) =0,
	\end{equation*}
	where $f_2 =div(H)- (m-n)\|H\|^2 + \mu - g(H, X)$ is a smooth function on $M$, which implies $(ii)$. This completes the proof.
\end{proof}
\begin{theorem}
	Let $(M, g, \xi, \mu)$ be a Ricci soliton with the potential vector field $\xi \in \Gamma(TM)$ and $F:(M^m, g) \to (N^n, h)$ be a totally geodesic horizontally conformal submersion between Riemannian manifolds. Then following statements are true:
	\begin{enumerate}[(i)]
		\item If $\xi= Z \in \Gamma(KerF_\ast)^\bot$, then $N$ is a Ricci soliton.
		\item If $\xi= U \in \Gamma(KerF_\ast)$, then $N$ is an Einstein.
	\end{enumerate}
\end{theorem}
\begin{proof}
	Since $(M, g, \xi, \mu)$ is a Ricci soliton then by (\ref{eqn1.2}), we have
	\begin{equation*}
		\frac{1}{2} (L_\xi g)(X,Y) + Ric(X,Y) + \mu g(X,Y) = 0,~\text{for}~ X, Y \in \Gamma(KerF_\ast)^\bot,
	\end{equation*}
	which can be written as
	\begin{equation}\label{eqn4.7}
		\frac{1}{2} \{g(\nabla_X \xi, Y) + g(\nabla_Y \xi, X)\} + Ric(X,Y) + \mu g(X,Y) = 0.
	\end{equation}
	If $\xi= Z \in \Gamma(KerF_\ast)^\bot$, then using (\ref{eqn2.1}) in (\ref{eqn4.7}), we get
	\begin{equation}\label{eqn4.8}
		\begin{array}{ll}
			\frac{1}{2\lambda^2} \{h(F_\ast (\nabla_X Z), F_\ast Y) + h(F_\ast(\nabla_Y Z), F_\ast X)\} \\+ Ric(X,Y) + \frac{\mu}{\lambda^2} h(F_\ast X,F_\ast Y) = 0.
		\end{array}
	\end{equation}
	Using (\ref{eqn2.10}) in (\ref{eqn4.8}), we get
	\begin{equation*}
		\begin{array}{ll}
			\frac{1}{2\lambda^2} \{ h(\nabla_{\tilde{X}}^N \tilde{Z}-( \nabla F_\ast)(X, Z), \tilde{Y}) \\+ h(\nabla_{\tilde{Y}}^N \tilde{Z}-( \nabla F_\ast)(Y, Z), \tilde{X})\} + Ric(X,Y) + \frac{\mu}{\lambda^2} h(\tilde{X}, \tilde{Y}) = 0,
		\end{array}
	\end{equation*}
	where $X$, $Y$ and $Z$ are the horizontal lift of $\tilde{X}$, $\tilde{Y}$ and $\tilde{Z}$, respectively. Since $F$ is totally geodesic, second fundamental form is zero and using Corollary \ref{cor3.3} in above equation, we get
	\begin{equation*}
		\begin{array}{ll}
			\frac{1}{2\lambda^2} \{ h(\nabla_{\tilde{X}}^N \tilde{Z}, \tilde{Y}) + h(\nabla_{\tilde{Y}}^N \tilde{Z}, \tilde{X})\} +\frac{1}{\lambda^2} Ric^N(\tilde{X},\tilde{Y}) + \frac{\mu}{\lambda^2} h(\tilde{X},\tilde{Y}) = 0,
		\end{array}
	\end{equation*}
	which implies $(i)$.\\Now, if $\xi= U \in \Gamma(KerF_\ast),$ then (\ref{eqn4.7}) implies 
	\begin{equation*}
		-\frac{1}{2} \{g(\nabla_X Y, U) + g(\nabla_Y X, U)\} + Ric(X,Y) + \mu g(X,Y) = 0.
	\end{equation*}
	Using (\ref{eqn2.1}) and (\ref{eqn2.6}) in above equation, we get
	\begin{equation}\label{eqn4.9}
		-\frac{1}{2} g(A_X Y + A_Y X, U) + Ric(X,Y) + \mu\frac{1}{\lambda^2} h(\tilde{X},\tilde{Y}) = 0,
	\end{equation}
	where $X$ and $Y$ are the horizontal lift of $\tilde{X}$ and $\tilde{Y}$, respectively. Since $F$ is totally geodesic then using Corollary \ref{cor3.3} in (\ref{eqn4.9}), we get
	\begin{equation*}
		Ric^N(\tilde{X},\tilde{Y}) + \mu h(\tilde{X},\tilde{Y}) = 0,
	\end{equation*}
	which implies $(ii)$. This completes the proof.
\end{proof}
\begin{theorem}
	Let $(M, g, \xi, \mu)$ be a Ricci soliton with the potential vector field $\xi \in \Gamma(TM)$ and $F:(M^m, g) \to (N^n, h)$ be a horizontally homothetic conformal submersion between Riemannian manifolds such that fibers of $F$ are totally geodesic and $(KerF_\ast)^\bot$ is integrable. Then following statements are true:
	\begin{enumerate}[(i)]
		\item If $\xi= Z \in \Gamma(KerF_\ast)^\bot$, then $N$ is an almost Ricci soliton.
		\item If $\xi= U \in \Gamma(KerF_\ast)$, then $N$ is an Einstein.
	\end{enumerate}
\end{theorem}
\begin{proof}
	If $\xi= Z \in \Gamma(KerF_\ast)^\bot$ then by using Lemma \ref{lem2.1} in (\ref{eqn4.8}) and condition of homothetic map, we get
	\begin{equation}\label{eqn4.10}
		\begin{array}{ll}
			\frac{1}{2\lambda^2} \{ h(\nabla_{\tilde{X}}^N \tilde{Z}, \tilde{Y}) + h(\nabla_{\tilde{Y}}^N \tilde{Z}, \tilde{X})\} + Ric(X,Y) + \frac{\mu}{\lambda^2} h(\tilde{X},\tilde{Y}) = 0.
		\end{array}
	\end{equation}
	Using Corollary \ref{cor3.2} in (\ref{eqn4.10}), we get	
	\begin{equation}\label{eqn4.11}
		\begin{array}{ll}
			\frac{1}{2\lambda^2} \{ h(\nabla_{\tilde{X}}^N \tilde{Z}, \tilde{Y}) + h(\nabla_{\tilde{Y}}^N \tilde{Z}, \tilde{X})\} + g(X, Y) div(H')+ \frac{1}{\lambda^2} Ric^N (\tilde{X}, \tilde{Y}) \\- \frac{1}{4} \lambda^4 g(X, Y) |\nabla_{\nu} \frac{1}{\lambda^2}|^2 + \frac{\mu}{\lambda^2} h(\tilde{X},\tilde{Y}) + \frac{n \lambda^2}{2} g(X, Y) \left(H' \frac{1}{\lambda^2}\right) = 0.
		\end{array}
	\end{equation}
	Using (\ref{eqn2.1}) in (\ref{eqn4.11}), we get
	\begin{equation*}
		\begin{array}{ll}
			\frac{1}{2} \{ h(\nabla_{\tilde{X}}^N \tilde{Z}, \tilde{Y}) + h(\nabla_{\tilde{Y}}^N \tilde{Z}, \tilde{X})\} + Ric^N (\tilde{X}, \tilde{Y}) + f_3 h(\tilde{X},\tilde{Y}) = 0,
		\end{array}
	\end{equation*}
	where $f_3 = \mu 
	+ div(H') - \frac{1}{4} \lambda^4 |\nabla_{\nu} \frac{1}{\lambda^2}|^2 + \frac{n \lambda^2}{2} \left(H' \frac{1}{\lambda^2}\right)$ is a smooth function on $M$, which implies $(i)$.\\ Now, if $\xi= U \in \Gamma(KerF_\ast),$ then using (\ref{eqn4.9}) and Corollary \ref{cor3.2}, we get
	\begin{equation*}
		\begin{array}{ll}
			-\frac{1}{2} g(A_X Y + A_Y X, U) + g(X, Y) div(H')+ \frac{1}{\lambda^2} Ric^N (\tilde{X}, \tilde{Y}) \\- \frac{1}{4} \lambda^4 g(X, Y) |\nabla_{\nu} \frac{1}{\lambda^2}|^2 + \mu\frac{1}{\lambda^2} h(\tilde{X},\tilde{Y}) + \frac{n \lambda^2}{2} g(X, Y) \left(H' \frac{1}{\lambda^2}\right) = 0.
		\end{array}
	\end{equation*}
	Using (\ref{eqn3.3}) and (\ref{eqn2.1}) in above equation, we get
	\begin{equation*}
		Ric^N (\tilde{X}, \tilde{Y})+ f_4 h(\tilde{X}, \tilde{Y})= 0,
	\end{equation*}
	where $f_4 =\frac{\lambda^2}{2} g(\nabla_\nu \frac{1}{\lambda^2}, U) + div(H')- \frac{1}{4} \lambda^4 |\nabla_{\nu} \frac{1}{\lambda^2}|^2 + \mu + \frac{n \lambda^2}{2} \left(H' \frac{1}{\lambda^2}\right)$ is a smooth function on $M$, which implies $(ii)$. This completes the proof.
\end{proof}
\begin{theorem}
	Let $(M, g, Z, \mu)$ be a Ricci soliton with the potential vector field $Z \in \Gamma(KerF_\ast)^\bot$ and $F:(M^m, g) \to (N^n, h)$ be a horizontally homothetic conformal submersion from a Riemannian manifold $M$ to an Einstein manifold $N$ such that fibers of $F$ are totally geodesic and $(KerF_\ast)^\bot$ is integrable. Then $Z$ is conformal vector field on $(KerF_\ast)^\bot$.
\end{theorem}
\begin{proof}
	Since $(M, g, Z, \mu)$ is a Ricci soliton then by (\ref{eqn1.2}), we have
	\begin{equation*}
		\frac{1}{2} (L_Z g)(X,Y) + Ric(X,Y) + \mu g(X,Y) = 0,~\text{for}~ X, Y, Z \in \Gamma(KerF_\ast)^\bot.
	\end{equation*}
	Using Corollary \ref{cor3.2} in above equation, we get
	\begin{equation}
		\begin{array}{ll}\label{eqn4.12}
			\frac{1}{2} (L_Z g)(X,Y) + g(X, Y) div(H')+ \frac{1}{\lambda^2} Ric^N (\tilde{X}, \tilde{Y})\\- \frac{1}{4} \lambda^4 g(X, Y) |\nabla_{\nu} \frac{1}{\lambda^2}|^2 + \mu g(X,Y) + \frac{n \lambda^2}{2} g(X, Y) \left(H' \frac{1}{\lambda^2}\right) = 0.
		\end{array}
	\end{equation}
	Using (\ref{eqn2.1}) in (\ref{eqn4.12}), we get
	\begin{equation}\label{eqn4.13}
		\frac{1}{2} (L_Z g)(X,Y) + \frac{1}{\lambda^2}\{ Ric^N (\tilde{X}, \tilde{Y}) + f_3 h(\tilde{X},\tilde{Y}) \}= 0,
	\end{equation}
	where $f_3 = div(H') - \frac{1}{4} \lambda^4 |\nabla_{\nu} \frac{1}{\lambda^2}|^2 + \mu + \frac{n \lambda^2}{2} \left(H' \frac{1}{\lambda^2}\right)$ is a smooth function on $M$. Since $N$ is an Einstein manifold, putting $Ric^N (\tilde{X}, \tilde{Y}) = f_3 h(\tilde{X},\tilde{Y})$ in (\ref{eqn4.13}), we get
	\begin{equation*}
		\frac{1}{2} (L_Z g)(X,Y) + \frac{2}{\lambda^2}f_3 h(\tilde{X},\tilde{Y}) = 0.
	\end{equation*}
	Using (\ref{eqn2.1}) in above equation, we get
	\begin{equation*}
		\frac{1}{2} (L_Z g)(X,Y) + 2f_3 g(X,Y) = 0,
	\end{equation*}
	which completes the proof.
\end{proof}
\begin{theorem}
	Let $(M, g, Z, \mu)$ be a Ricci soliton with the potential vector field $Z \in \Gamma(KerF_\ast)^\bot$ and $F:(M^m, g) \to (N^n, h)$ be a totally geodesic horizontally conformal submersion from a Riemannian manifold $M$ to a connected Einstein manifold $N$. Then $\tilde{Z}$ is conformal vector field on $N$.
\end{theorem}
\begin{proof}
	Since $(M, g, Z, \mu)$ is a Ricci soliton then by (\ref{eqn1.2}), we have
	\begin{equation*}
		\frac{1}{2} (L_Z g)(X,Y) + Ric(X,Y) + \mu g(X,Y) = 0,~\text{for}~ X, Y, Z \in \Gamma(KerF_\ast)^\bot.
	\end{equation*}
	Using (\ref{eqn2.1}) and Corollary \ref{cor3.3} in above equation, we get
	\begin{equation}\label{eqn4.14}
		\begin{array}{ll}
			\frac{1}{2} (L_Z g)(X,Y) +\frac{1}{\lambda^2} Ric^N (\tilde{X}, \tilde{Y}) + \frac{\mu}{\lambda^2} h(\tilde{X}, \tilde{Y}) = 0.
		\end{array}
	\end{equation}
	Since $N$ is an Einstein manifold, putting $Ric^N (\tilde{X}, \tilde{Y}) = \mu h(\tilde{X},\tilde{Y})$ in (\ref{eqn4.14}), we get
	\begin{equation*}
		\frac{1}{2} (L_Z g)(X,Y) + \frac{2 \mu}{\lambda^2} h(\tilde{X},\tilde{Y}) = 0,
	\end{equation*}
	which implies
	\begin{equation}\label{eqn4.15}
		\frac{1}{2} \{g(\nabla_X Z, Y) + g(\nabla_Y Z ,X)\} + \frac{2 \mu}{\lambda^2} h(\tilde{X},\tilde{Y}) = 0.
	\end{equation}
	Using (\ref{eqn2.1}) in (\ref{eqn4.15}), we get
	\begin{equation*}
		\frac{1}{2\lambda^2} \{h(F_\ast (\nabla_X Z), F_\ast Y) + h(F_\ast(\nabla_Y Z), F_\ast X)\} + \frac{2 \mu}{\lambda^2} h(\tilde{X},\tilde{Y}) = 0.
	\end{equation*}
	Since $F$ is totally geodesic and using (\ref{eqn2.10}) in above equation, we get
	\begin{equation*}
		\frac{1}{2} \{ h(\nabla_{\tilde{X}}^N \tilde{Z}, \tilde{Y}) + h(\nabla_{\tilde{Y}}^N \tilde{Z}, \tilde{X})\} + 2 \mu  h(\tilde{X},\tilde{Y}) = 0,
	\end{equation*}
	which implies
	\begin{equation*}
		\frac{1}{2} (L_{\tilde{Z}} h)(\tilde{X},\tilde{Y}) + 2 \mu  h(\tilde{X},\tilde{Y}) = 0.
	\end{equation*}
	Thus the proof is completed.
\end{proof}
\begin{theorem}
	Let $(M, g, \xi, \mu)$ be a Ricci soliton with the potential vector field $\xi \in \Gamma(TM)$ and $F:(M^m, g) \to (N^n, h)$ be a totally geodesic horizontally conformal submersion between Riemannian manifolds. Then $M$ has constant scalar curvature given by $-\mu m$.
\end{theorem}
\begin{proof}  Since $(M,g,\xi,\mu)$ is a Ricci soliton then from (\ref{eqn1.2}), we get
	\begin{equation}\label{eqn4.16} 
		\frac{1}{2} \{ g (\nabla_X \xi, Y) + g (\nabla_Y \xi, X) \} + Ric (X, Y)+ \mu g(X,Y) = 0,
	\end{equation}
	for any $X,Y,\xi \in \Gamma(TM)$. Now, we decompose $X,Y$ and $\xi$ such that $X= \nu X + \mathcal{H} X$, $Y= \nu Y + \mathcal{H} Y $ and $\xi= \nu \xi + \mathcal{H}  \xi $, where $\nu$  and $\mathcal{H} $ denote for component of $KerF_\ast$ and $(KerF_\ast)^\bot$, respectively. Then (\ref{eqn4.16}) can be written as
	\begin{eqnarray}\label{eqn4.17} 
		\begin{array}{ll}
			\frac{1}{2} \{ g (\nabla_{\nu X + \mathcal{H} X } \nu \xi + \mathcal{H} \xi, \nu Y + \mathcal{H} Y) + g (\nabla_{\nu Y + \mathcal{H} Y} \nu \xi + \mathcal{H} \xi, \nu X + \mathcal{H} X) \} \\+ Ric (\nu X, \nu Y)+ Ric (\mathcal{H} X, \mathcal{H} Y) + Ric (\nu X, \mathcal{H} X)+ Ric (\nu Y, \mathcal{H} Y) \\+  \mu \{ g (\nu X, \nu Y)+ g (\mathcal{H} X, \mathcal{H} Y)\} = 0.
		\end{array}
	\end{eqnarray}
	Taking trace of (\ref{eqn4.17}), we get
	\begin{eqnarray*} 
		\begin{array}{ll}
			\frac{1}{2} \Big \{ 2 \sum\limits_{i=n+1}^{m} g (\nabla_{U_i} U_i, U_i) + 2 \sum\limits_{j=1}^{n} g (\nabla_{X_j} X_j, X_j) \Big \} +\sum\limits_{i=n+1}^{m}Ric (U_i, U_i) \\+ \sum\limits_{j=1}^{n}  Ric (X_j, X_j) + 2 \underset{i, j}{\overset{}{\sum}} Ric (U_i, X_j)+  \mu \Big \{ \sum\limits_{i=n+1}^{m}g (U_i, U_i)+ \sum\limits_{j=1}^{n}  g (X_j, X_j)\Big \} = 0,
		\end{array}
	\end{eqnarray*}
	where $\{U_i\}_{n+1 \leq i \leq m}$ and $\{X_j\}_{1 \leq j \leq n}$ are orthonormal bases of $KerF_\ast$ and $(KerF_\ast)^\bot$, respectively. Since $F$ is totally geodesic then using Corollary \ref{cor3.3}, in above equation, we get
	\begin{equation}\label{eqn4.18} 
		\begin{array}{ll}
			\sum\limits_{i=n+1}^{m}g (\nabla_{U_i} U_i, U_i) + \sum\limits_{j=1}^{n} g (\nabla_{X_j} X_j, X_j) \\+ s^{KerF_\ast} + \frac{1}{\lambda^2}s^{N} + \mu (m-n+n) = 0,
		\end{array}
	\end{equation}
	where $s^{KerF_\ast}$ and $s^N$ denote the scalar curvatures of $KerF_\ast$ and $N$, respectively. \\Since $\nabla$ is metric connection on $M$ then from (\ref{eqn4.18}), we get
	\begin{equation*} 
		\frac{1}{2} \Big \{ \sum\limits_{i=n+1}^{m}\nabla_{U_i} (g (U_i, U_i)) + \sum\limits_{j=1}^{n} \nabla_{X_j} (g ( X_j, X_j)) \Big \} + s^{KerF_\ast}+ \frac{1}{\lambda^2}s^{N} + \mu m = 0.
	\end{equation*}
	Using Theorem \ref{thm3.4} in above equation, we get
	\begin{equation*} 
		\frac{1}{2} \Big \{ \sum\limits_{i=n+1}^{m}{U_i}(g(U_i, U_i)) + \sum\limits_{j=1}^{n} {X_j} (g(X_j,X_j)) \Big \} + s + \mu m = 0.
	\end{equation*}
	Thus $s + \mu m = 0$, where $s$ is the scalar curvature on $M$. Hence $s= -\mu m$, which completes the proof.
\end{proof}
\begin{definition} (\cite{Sahin_2017}, p. 45) 
	Let $F:(M^m,g) \to (N^n,h)$ be a smooth map between Riemannian manifolds. Then $F$ is harmonic if and only if the tension field $\tau(F)$ of $F$ vanishes at each point $p\in M$.
\end{definition}
\begin{theorem}
	Let $(M, g, W, \mu)$ be a Ricci soliton with the potential vector field $W \in \Gamma(KerF_\ast)$ and $F:(M^m, g) \to (N^n, h)$ be a horizontally homothetic conformal submersion between Riemannian manifolds such that fibers of $F$ are totally umbilical and $(KerF_\ast)^\bot$ is totally geodesic. Then $F$ is harmonic if and only if the scalar curvature of $KerF_\ast$ is $-\mu(m-n)$.
\end{theorem}
\begin{proof}
	Let $F:(M^m,g) \to (N^n,h)$ be a conformal submersion between Riemannian manifolds. Then we have \cite{Meena_2022}
	\begin{equation}\label{eqn4.19}
		\tau(F) = (n-2) \frac{\lambda^2}{2} F_\ast\left(\nabla_{\mathcal{H}}\frac{1}{\lambda^2} \right) -(m-n) F_\ast (H).
	\end{equation}
	Now putting $\xi =W$ in (\ref{eqn4.6}) and then taking trace, we get
	\begin{equation*}
		\begin{array}{ll}
			\frac{1}{2} \{2 \sum\limits_{i=n+1}^{m} g(\nabla_{U_i} U_i , U_i)\} + \sum\limits_{i=n+1}^{m} Ric^\nu (U_i, U_i) \\- (m-n)\|H\|^2\sum\limits_{i=n+1}^{m}g(U_i, U_i) + div(H) \sum\limits_{i=n+1}^{m}g(U_i, U_i) + \mu \sum\limits_{i=n+1}^{m}g(U_i, U_i) = 0,
		\end{array}
	\end{equation*}
	where $\{U_i\}_{n+1 \leq i \leq m}$ be an orthonormal basis of $KerF_\ast$, which implies
	\begin{equation}\label{eqn4.20}
		s^{KerF_\ast} + (m-n) \mu - (m-n)^2 \|H\|^2 + (m-n) div(H)= 0,
	\end{equation}
	where $s^{KerF_\ast}$ is the scalar curvature of $KerF_\ast$. Since $s^{KerF_\ast} = -\mu (m-n)$ then (\ref{eqn4.20}) implies
	\begin{equation*}
		div (H)-(m-n) \|H\|^2 = 0 \iff H=0.
	\end{equation*}
	Since $F$ is homothetic conformal submersion, then by (\ref{eqn4.19}), we get $\tau(F) =0 \iff H=0$. This completes the proof.
\end{proof}
\section{Examples}\label{sec5}
In this section, we give four non-trivial examples to support the theory of the paper.

First, we give an example of a homothetic horizontally conformal submersion $F$ such that fibers of $F$ are totally geodesic and $(KerF_\ast)^\bot$ is integrable.
\begin{ex}\label{example5.1}
	Let $M= \{ (x_1, x_2) \in \mathbb{R}^{2} : x_1 \neq 0, x_2 \neq 0\}$ be a Riemannian manifold with Riemannian metric $g$ on $M$ given by $g= e^{-2x_2}dx_1^2 + dx_2^2$. Let $N=\{ y_1 \in \mathbb{R} : y_1 \neq 0\}$ be a Riemannian manifold with Riemannian metric $h$ on $N$ given by $h= dy_1^2$. Consider a map $F : (M, g) \rightarrow (N, h) $ defined by
	\begin{equation*}
		F(x_1,x_2)= x_1.
	\end{equation*}
	By direct computations
	\begin{equation*}	
		KerF_\ast = Span \{ U = e_2\},
	\end{equation*}
	and 
	\begin{equation*}
		(KerF_\ast)^\bot = Span \{ X = e_1\},
	\end{equation*}
	where $\left\{ e_1 = \frac{\partial}{\partial x_1}, e_2 = \frac{\partial}{\partial x_2}\right\}$, $\Big\{ e_1' = \frac{\partial}{\partial y_1} \Big\}$ are bases of  $T_pM$ and $T_{F (p)}N$ respectively, for any $ p\in M$. In addition, we can see that $F_\ast (X = e_1) = e_1'$, $F_\ast (U = e_2) = 0$ and $\lambda^2 g(X,X)= h(F_\ast X, F_\ast X)$ for $\lambda = e^{x_2}$ and $X \in \Gamma(KerF_\ast)^\bot.$ Thus $F$ is a horizontally conformal submersion with $RangeF_\ast= Span \{e_1' \}$. Now, we compute the Christoffel symbols for the metric $g$
	\begin{equation*}
		\Gamma_{11}^{1} = 0, \Gamma_{11}^2 = e^{-2x_2}, \Gamma_{12}^1 = -1 = \Gamma_{21}^1, \Gamma_{12}^2 = 0= \Gamma_{21}^2, \Gamma_{22}^1 = 0, \Gamma_{22}^2 = 0.
	\end{equation*}
	Then using above equation, we get
	\begin{equation}\label{eqn5.1}
		\begin{array}{ll}
			\nabla_{e_1} e_1 = \nabla_{X} X = e^{-2x_2} e_2,\\ \nabla_{e_2} e_2 = \nabla_{U} U = 0,\\ \nabla_{e_1} e_2 = \nabla_{X} U = -e_1, \\ \nabla_{e_2} e_1 = \nabla_{U} X = -e_1.
		\end{array}
	\end{equation}
	Using (\ref{eqn2.4}), (\ref{eqn2.5}), (\ref{eqn2.6}) and (\ref{eqn5.1}), we get
	\begin{equation*}
		A_X X = e^{-2x_2} e_2,~ T_U U =0,~ [X, X] = 0, ~\text{and}~ \mathcal{H}[X, U] = 0.
	\end{equation*}
	Also, 
	\begin{equation*}
		\mathcal{H}(grad~ \lambda) = \mathcal{H}(grad~ e^{x_2}) = \mathcal{H}(e^{x_2} e_2) =0.
	\end{equation*}
	Thus, $F$ is homothetic horizontally conformal submersion such that fibers of $F$ are totally geodesic and $(KerF_\ast)^\bot$ is integrable. 
	Now, we will show that total manifold $M$ admits a Ricci soliton, i.e.
	\begin{equation}\label{eqn6.1} 
		\frac{1}{2} (L_{Z_1} g)(X_1,Y_1) + Ric (X_1,Y_1) + \mu g(X_1,Y_1) = 0,
	\end{equation}
	for any $X_1, Y_1, Z_1 \in \Gamma(TM)$. We know that
	\begin{equation}\label{eqn6.2} 
		\frac{1}{2} (L_{Z_1} g)(X_1,Y_1) = \frac{1}{2} \Big \{ g(\nabla_{X_1} Z_1,Y_1) + g(\nabla_{Y_1} Z_1,X_1) \Big \}.
	\end{equation}
	Since dimension of $KerF_\ast$ and $(KerF_\ast)^\bot$ is one, we can decompose $X_1, Y_1$ and $Z_1$ such that $X_1= a_1e_1+a_2e_2$, $Y_1=a_3e_1+a_4e_2$ and $Z_1= a_5e_1+a_6e_2$, where
	$\{a_i\}_{1 \leq i \leq 6} \in \mathbb{R}$ are some scalars. Then from (\ref{eqn6.2}), we get
	\begin{equation*}
		\begin{array}{ll}
			\frac{1}{2} (L_{Z_1} g)(X_1,Y_1)=& \frac{1}{2} \Big \{ g(\nabla_{a_1e_1+a_2e_2} a_5e_1+a_6e_2,a_3e_1+a_4e_2) \\&+ g(\nabla_{a_3e_1+a_4e_2} a_5e_1+a_6e_2,a_1e_1+a_2e_2) \Big \}.
		\end{array}
	\end{equation*}
	Since $\nabla$ is linear connection and using (\ref{eqn5.1}) in above equation, we get
	\begin{equation}\label{eqn6.5}
		\frac{1}{2} (L_{Z_1} g)(X_1,Y_1)= \frac{1}{2}\{e^{-2x_2}a_1a_4a_5+ e^{-2x_2}a_2a_4a_5 -2e^{-2x_2}a_1a_3a_6 \}.
	\end{equation}
	Also,
	\begin{equation}\label{eqn6.6} 
		g(X_1,Y_1)= g(a_1e_1+a_2e_2,a_3e_1+a_4e_2) = (a_1a_3e^{-2x_2} + a_2a_4),
	\end{equation}
	and
	\begin{equation*}
		Ric(X_1,Y_1) = Ric(a_1e_1+a_2e_2,a_3e_1+a_4e_2),
	\end{equation*}
	which implies
	\begin{equation}\label{eqn6.7} 
		Ric(X_1,Y_1)= a_1a_3Ric(e_1,e_1) + (a_1a_4+a_2a_3) Ric(e_1,e_2) +a_2a_4 Ric(e_2,e_2).
	\end{equation}
	Using Proposition \ref{prop3.2} or Corollary \ref{Cor3.1}, we get
	\begin{equation}\label{eqn6.8} 
		Ric(e_1, e_1)= -e^{-6x_2}-e^{-2x_2} -e^{-4x_2} + 1,
	\end{equation}
	\begin{equation}\label{eqn6.9} 
		Ric(e_2, e_2)=-2e^{-2 x_2} - 1,
	\end{equation}
	and
	\begin{align}\label{eqn6.10} 
		Ric(e_1, e_2)=0.
	\end{align}
	Using (\ref{eqn6.8}), (\ref{eqn6.9}) and (\ref{eqn6.10}) in (\ref{eqn6.7}), we get
	\begin{equation}\label{eqn6.11} 
		Ric(X_1,Y_1)= -a_1a_3e^{-6x_2}-a_1a_3e^{-2x_2} -a_1a_3e^{-4x_2} + a_1a_3- 2a_2 a_4 e^{-2 x_2} - a_2a_4.
	\end{equation}
	Now, using (\ref{eqn6.5}), (\ref{eqn6.6}) and (\ref{eqn6.11}) in (\ref{eqn6.1}), we obtain that $M$ admits a Ricci soliton for
	{\tiny\begin{equation*}
		\mu =  \frac{\{2a_1a_3 (a_6 + e^{-2x_2}) -a_4a_5 (a_1 + a_2)+2a_1a_3 (1+e^{-4x_2})\}e^{-2x_2} -2a_1a_3 + 2a_2a_4 (1+ 2e^{-2x_2} )}{2(a_1a_3e^{-2x_2} + a_2a_4)},
	\end{equation*}
	}where $a_1a_3e^{-2x_2} \neq - a_2a_4$. Since all $a_i \in \mathbb{R}$, for some choices of $a_i$s Ricci soliton will be shrinking, expanding or steady according to $\mu < 0$, $\mu >0$ or $\mu = 0$.	
\end{ex}
Further, we give an example of horizontally conformal submersion $F$ such that fibers of $F$ are totally geodesic and totally umbilical. In addition, $(KerF_\ast)^\bot$ is integrable and totally geodesic.
\begin{ex}
	Let $M= \{ (x_1, x_2) \in \mathbb{R}^{2} : x_1 \neq 0\}$ be a Riemannian manifold with Riemannian metric $g$ on $M$ given by $g= e^{2x_1}dx_1^2 + dx_2^2$. Let $N=\{ y_1 \in \mathbb{R} : y_1 \neq 0\}$ be a Riemannian manifold with Riemannian metric $h$ on $N$ given by $h= dy_1^2$. Consider a map $F : (M, g) \rightarrow (N, h) $ defined by
	\begin{equation*}
		F(x_1,x_2)= x_1.
	\end{equation*}
	By direct computations
	\begin{equation*}	
		KerF_\ast = Span \{ U = e_2\},
	\end{equation*}
	and 
	\begin{equation*}
		(KerF_\ast)^\bot = Span \{ X = e_1\},
	\end{equation*}
	where $\left\{ e_1 = \frac{\partial}{\partial x_1}, e_2 = \frac{\partial}{\partial x_2}\right\}$, $\Big\{ e_1' = \frac{\partial}{\partial y_1}\Big\}$ are bases of  $T_pM$ and $T_{F (p)}N$ respectively, for any $ p\in M$. In addition, we can see that $F_\ast (X = e_1) = e_1'$, $F_\ast (U = e_2) = 0$ and $\lambda^2 g(X,X)= h(F_\ast X, F_\ast X)$ for $\lambda = e^{-x_1}$ and $X \in \Gamma(KerF_\ast)^\bot.$ Thus $F$ is a horizontally conformal submersion with $RangeF_\ast= Span \{e_1' \}$. Now, we compute the Christoffel symbols for the metric $g$
	\begin{equation*}
		\Gamma_{11}^{1} = 1, \Gamma_{11}^2 = 0, \Gamma_{12}^1 = 0= \Gamma_{21}^1, \Gamma_{12}^2 = 0= \Gamma_{21}^2, \Gamma_{22}^1 = 0, \Gamma_{22}^2 = 0.
	\end{equation*}
	Then using above equation, we get
	\begin{equation}\label{eqn5.13}
		\begin{array}{ll}
			\nabla_{e_1} e_1 = \nabla_{X} X = e_1,\\ \nabla_{e_2} e_2 = \nabla_{U} U = 0,\\ \nabla_{e_1} e_2 = \nabla_{X} U = 0, \\ \nabla_{e_2} e_1 = \nabla_{U} X = 0.
		\end{array}
	\end{equation}
	Using (\ref{eqn2.4}), (\ref{eqn2.5}), (\ref{eqn2.6}) and (\ref{eqn5.13}), we get
	\begin{equation*}
		A_X X = 0,~ T_U U =0~\text{and}~ [X, X] = 0.
	\end{equation*}
	Also, 
	\begin{equation*}
		\mathcal{H}(grad~ \lambda) \neq 0~\text{and}~\nu (grad~ \lambda) = 0.
	\end{equation*}
	Thus, $F$ is horizontally conformal submersion such that fibers of $F$ are totally geodesic and totally umbilical also. In addition $(KerF_\ast)^\bot$ is integrable and totally geodesic.
	
	Now, we will show that total manifold $M$ admits a Ricci soliton, i.e. (\ref{eqn6.1}) holds.
	Since dimension of $KerF_\ast$ and $(KerF_\ast)^\bot$ is one, we can decompose $X_1, Y_1$ and $Z_1$ such that $X_1= a_1e_1+a_2e_2$, $Y_1=a_3e_1+a_4e_2$ and $Z_1= a_5e_1+a_6e_2$, where
	$\{a_i\}_{1 \leq i \leq 6} \in \mathbb{R}$ are some scalars. Then using (\ref{eqn6.2}) and (\ref{eqn5.13}), we get
	\begin{equation}\label{eqn5.14} 
		\frac{1}{2} (L_{Z_1} g)(X_1,Y_1)= a_1a_3a_5 e^{2 x_1}.
	\end{equation}
	Also,
	\begin{equation}\label{eqn5.15} 
		g(X_1,Y_1)= g(a_1e_1+a_2e_2,a_3e_1+a_4e_2) = (a_1a_3e^{2x_1} + a_2a_4),
	\end{equation}
	and by using Proposition \ref{prop3.2} in (\ref{eqn6.7}), we get
	\begin{equation}\label{eqn5.17} 
		Ric(X_1,Y_1)= a_1 a_3 (1-2e^{2 x_1})(e^{2 x_1} - 1).
	\end{equation}
	Now, using (\ref{eqn5.14}), (\ref{eqn5.15}) and (\ref{eqn5.17}) in (\ref{eqn6.1}), we obtain $$\mu =  \frac{a_1 a_3 (1-2e^{2 x_1})(1 - e^{2 x_1})-a_1a_3a_5e^{2x_1}}{a_1a_3e^{2x_1} + a_2a_4}$$ where $a_1a_3e^{2x_1} \neq - a_2a_4$. Since all $a_i \in \mathbb{R}$, for some choices of $a_i$s Ricci soliton will be shrinking, expanding or steady according to $\mu < 0$, $\mu >0$ or $\mu = 0$.	
\end{ex}
Next, we give an example of a horizontally conformal submersion such that fibers of $F$ are totally umbilical, and $(KerF_\ast)^\bot$ is integrable and totally geodesic.
\begin{ex}\label{example5.2}
	Let $M= \{ (x_1, x_2, x_3) \in \mathbb{R}^{3} :x_3 > 1\}$ be a Riemannian manifold with Riemannian metric $g$ on $M$ given by $g= x_3^{-2}dx_1^2 + x_3^{-2}dx_2^2 +x_3^{-2}dx_3^2$. Let $N=\{ (y_1, y_2)\in \mathbb{R}^2: y_2 > 1\}$ be a Riemannian manifold with Riemannian metric $h$ on $N$ given by $h= dy_1^2 + dy_2^2$. Consider a map $F : (M, g) \rightarrow (N, h) $ defined by
	\begin{equation*}
		F(x_1,x_2, x_3)= (x_2, x_3).
	\end{equation*}
	By direct computations
	\begin{equation*}	
		KerF_\ast = Span \{ U = e_1\},
	\end{equation*}
	and 
	\begin{equation*}
		(KerF_\ast)^\bot = Span \{ X_1 = e_2, X_2 = e_3\},
	\end{equation*}
	where $\left\{ e_1 = \frac{\partial}{\partial x_1}, e_2 =  \frac{\partial}{\partial x_2}, e_3 = \frac{\partial}{\partial x_3}\right\}$, $\Big\{ e_1' = \frac{\partial}{\partial y_1}, e_2' = \frac{\partial}{\partial y_2} \Big\}$ are bases of $T_pM$ and $T_{F (p)}N$ respectively, for any $ p\in M$. In addition, we can see that $F_\ast(X_1) = e_1'$, $F_\ast(X_2) = e_2'$, $F_\ast(U) = 0$ and $\lambda^2 g(X_i,X_j)= h(F_\ast X_i, F_\ast X_j)$ for $\lambda = x_3$ and $X_i, X_j \in \Gamma(KerF_\ast)^\bot.$ Thus $F$ is a horizontally conformal submersion with $RangeF_\ast= Span \{e_1', e_2'\}$. Now, we compute the Christoffel symbols for the metric $g$
	\begin{equation*}
		\begin{array}{ll}
			\Gamma_{13}^{1} = \Gamma_{31}^{1} = -x_3^{-1}, \Gamma_{23}^{2} = \Gamma_{32}^{2} = -x_3^{-1}, \Gamma_{11}^3 = x_3^{-1}, \Gamma_{22}^3 = x_3^{-1}, \Gamma_{33}^3 = -x_3^{-1},
		\end{array}
	\end{equation*}
	and remaining are zero. Then using above equation, we get
	\begin{equation}\label{eqn5.11}
		\begin{array}{ll}
			\nabla_{e_1} e_1 = \nabla_{U} U = x_3^{-1}e_3,\\ \nabla_{e_2} e_2 = \nabla_{X_1} X_1 = x_3^{-1}e_3,\\ \nabla_{e_1} e_2 = \nabla_{U} X_1 = 0, \\ \nabla_{e_2} e_1 = \nabla_{X_1} U = 0,\\ \nabla_{e_1} e_3 = \nabla_{U} X_2 = -x_3^{-1}e_1, \\ \nabla_{e_3} e_1 = \nabla_{X_2} U = -x_3^{-1}e_1,\\ \nabla_{e_2} e_3 = \nabla_{X_1} X_2 = -x_3^{-1}e_2, \\ \nabla_{e_3} e_2 = \nabla_{X_2} X_1 =-x_3^{-1}e_2,\\ \nabla_{e_3} e_3 = \nabla_{X_2} X_2 = -x_3^{-1}e_3.
		\end{array}
	\end{equation}
	Using (\ref{eqn2.4}), (\ref{eqn2.5}), (\ref{eqn2.6}) and (\ref{eqn5.11}), we get
	\begin{equation*}
		T_U U =e_3,
	\end{equation*}
	\begin{equation*}
		\nabla_X X = \nabla_{aX_1 + bX_2} aX_1 + bX_2 = (a^2-b^2)x_3^{-1}e_3-2abx_3^{-1}e_2 \implies A_X X = 0,
	\end{equation*}
	for $X \in \Gamma(KerF_\ast)^\bot~\text{and}~ a,b \in \mathbb{R}$, and
	\begin{equation*}
		[X_1, X_1] = [X_2, X_2] = 0= [X_1, X_2].	
	\end{equation*}
	Also, we see that $g(U, U) H = g(U, U) X = g(U, U) (aX_1 + bX_2) = x_3^{-2}e_3$ for $a =0$ and $b=1$. 
	Thus, $F$ is horizontally conformal submersion such that fibers of $F$ are totally umbilical, and $(KerF_\ast)^\bot$ is integrable and totally geodesic. 
	Now, by similar computation to Example \ref{example5.1} we can easily show that $M$ admits a Ricci soliton for some $\mu$, so we are omitting it.
\end{ex}
Finally, we give an example of a homothetic horizontally totally geodesic conformal submersion such that $KerF_\ast$ and $(KerF_\ast)^\bot$ are totally geodesic.
\begin{ex}\label{example5.3}
	Let $M= \{ (x_1, x_2, x_3) \in \mathbb{R}^{3}\}$ be a Riemannian manifold with Riemannian metric $g$ on $M$ given by $g= 4dx_1^2 + 4dx_2^2 +4dx_3^2$. Let $N=\{ (y_1, y_2)\in \mathbb{R}^2\}$ be a Riemannian manifold with Riemannian metric $h$ on $N$ given by $h= dy_1^2 + dy_2^2$. Consider a map $F : (M, g) \rightarrow (N, h) $ defined by
	\begin{equation*}
		F(x_1,x_2, x_3)= (x_1, x_3).
	\end{equation*}
	By direct computations
	\begin{equation*}	
		KerF_\ast = Span \{ U = e_2\},
	\end{equation*}
	and 
	\begin{equation*}
		(KerF_\ast)^\bot = Span \{ X_1 = e_1, X_2 = e_3\},
	\end{equation*}
	where $\left\{ e_1 = \frac{\partial}{\partial x_1}, e_2 = \frac{\partial}{\partial x_2}, e_3 = \frac{\partial}{\partial x_3}\right\}$, $\Big\{ e_1' = \frac{\partial}{\partial y_1}, e_2' = \frac{\partial}{\partial y_2} \Big\}$ are bases of  $T_pM$ and $T_{F (p)}N$ respectively, for any $ p\in M$. In addition, we can see that $F_\ast(X_1) = e_1'$, $F_\ast(X_2) = e_2'$, $F_\ast(U) = 0$ and $\lambda^2 g(X_i,X_j)= h(F_\ast X_i, F_\ast X_j)$ for $\lambda = \frac{1}{2}$ and $X_i, X_j \in \Gamma(KerF_\ast)^\bot.$ Thus $F$ is a horizontally conformal submersion with $RangeF_\ast= Span \{e_1', e_2'\}$. 
	Then using the Christoffel symbols $\Gamma^{i}_{jk} = 0$ for the metric $g$, we get
	\begin{equation}\label{eqn5.12}
		\begin{array}{ll}
			\nabla_{e_i} e_j = 0~\text{for $1 \leq i , j \leq 3$}.
		\end{array}
	\end{equation}
	Using (\ref{eqn2.4}), (\ref{eqn2.5}), (\ref{eqn2.6}) and (\ref{eqn5.12}), we get
	\begin{equation*}
		T_U U =0,
	\end{equation*}
	and
	\begin{equation*}
		\nabla_X X = \nabla_{aX_1 + bX_2} aX_1 + bX_2 = 0\implies A_X X = 0,
	\end{equation*}
	for $X \in \Gamma(KerF_\ast)^\bot~\text{and}~ a,b \in \mathbb{R}$. Thus, $F$ is homothetic horizontally conformal submersion such that $KerF_\ast$ and $(KerF_\ast)^\bot$ are totally geodesic. Thus by Theorem \ref{thm3.3}, $F$ is totally geodesic. Now, by similar computation to Example \ref{example5.1} we can easily show that $M$ admits a Ricci soliton for some $\mu$, so we are omitting it.
\end{ex}

\section*{Acknowledgments}
The authors are sincerely grateful to the referee(s) for his/her comments to
improve the paper.

\section*{Author contributions}
Conceptualization, Methodology, Investigation, writing-original draft preparation (K. Meena), Review (K. Meena and A. Yadav). All the authors read and agreed to the manuscript.

\section*{Funding}
The following financial supports were provided to the first author Kiran Meena during the preparation of this manuscript. (i). Council of Scientific and Industrial Research-Human Resource Development Group (CSIR-HRDG), New Delhi, India [File No.: 09/013(0887)/2019-EMR-I] during the doctoral study. (ii). Harish-Chandra Research Institute, Prayagraj, India [Offer Letter No.: HRI/133/1436 Dated 29 November 2022] during the post-doctoral study.

\section*{Data Availability Statement}
This manuscript has no associated data.

\section*{Declarations}
	\textbf{Conflict of interest}
		None of the authors have any relevant financial or non-financial competing interests.\\
		
	\noindent \textbf{Ethical Approval}
		The submitted work is original and not submitted to more than one journal for simultaneous consideration.\\
		
	\noindent \textbf{Publisher's Note}
	Springer Nature remains neutral with regard to jurisdictional claims in published maps and institutional affiliations.\\
	
	\noindent Springer Nature or its licensor (e.g. a society or other partner) holds exclusive rights to this article under a publishing agreement with the author(s) or other	rightsholder(s); author self-archiving of the accepted manuscript version of	this article is solely governed by the terms of such publishing agreement and	applicable law.\\

\end{document}